\documentclass[12pt]{article}
\usepackage{a4,times,amsmath,amscd,amsfonts,amsthm,latexsym, amssymb,
  oldgerm}

\newtheorem{thm}{\qquad\sc Theorem}\newcommand{\bt}{\begin{thm}}
  \newcommand{\et}{\end{thm}} \newtheorem{lemma}[thm]{\qquad \sc Lemma}
\newcommand{\bl}{\begin{lemma}} \newcommand{\el}{\end{lemma}}
\newtheorem{prop}[thm]{\qquad\sc Proposition} \newcommand{\bp}{\begin{prop}}
  \newcommand{\ep}{\end{prop}}
\newtheorem{rmk}[thm]{\qquad\sc Remark} \newcommand{\brmk}{\begin{rmk}}
  \newcommand{\ermk}{\end{rmk}}
\newcommand{\pf}{\qquad {\it Proof:}\quad }
\newtheorem{cor}[thm]{\qquad\sc Corollary}\newcommand{\bc}{\begin{cor}}
\newcommand{\ec}{\end{cor}}
\newtheorem{ex}[thm]{\qquad\sc Example} \newcommand{\bex}{\begin{ex}}
  \newcommand{\eex}{\end{ex}}

\begin{document}

\begin{center}
  {\LARGE Effective multiplicity one on ${\rm GL}_n$ and narrow zero-free regions for Rankin-Selberg $L$-functions}
\end{center}

\medskip

\begin{center}
  {\large {\sc Farrell Brumley}}
\end{center}

\bigskip

\begin{abstract}
\noindent We establish zero-free regions tapering as an inverse power of the analytic conductor for Rankin-Selberg $L$-functions 
on ${\rm GL}_n\times {\rm GL}_{n'}$.  Such zero-free regions are equivalent to commensurate lower bounds on the edge of the critical 
strip, and in the case of $L(s,\pi\times\widetilde{\pi})$, on the residue at $s=1$.  As an application we show that a cuspidal 
automorphic representation on ${\rm GL}_n$ is determined by a finite number of its Dirichlet series coefficients, and that this 
number grows at most polynomially in the analytic conductor.
\end{abstract}

\bigskip\bigskip\bigskip

\qquad Let ${\mathbb{A}}$ be the ring of adeles over a number field $F$ and let $\pi$ and $\pi '$ be two cuspidal representations 
of ${\rm GL}_n({\mathbb{A}})$ with restricted tensor product decompositions $\pi=\otimes_v \pi_v$ and $\pi '=\otimes_v \pi_v'$ over 
all places $v$ of $F$.  The strong multiplicity one theorem asserts that if $\pi_v\simeq\pi_v'$ for all but finitely many places 
$v$, then $\pi= \pi'$.  This was proven by Piatetski-Shapiro [PS] using the uniqueness of the Kirillov model and then by Jacquet 
and Shalika [J-S] using Rankin-Selberg $L$-functions.  Much more can be said however about the extent to which agreement of local 
factors on a suitable subset of the primes determines global equality.  For instance, Moreno has shown [Mo1] that for some finite 
$Y(\pi,\pi')$ the condition that $\pi_\textfrak{p}\simeq\pi_\textfrak{p}'$ for spherical non-archimedean $\textfrak{p}$ with 
absolute norm ${\rm N}\textfrak{p}\leq Y(\pi,\pi')$ is sufficient to imply $\pi=\pi'$.  From the analytic perspective, the 
crucial issue are the zeros of Rankin-Selberg $L$-functions: under GRH for both $L(s,\pi\times\tilde{\pi}')$ and 
$L(s,\pi\times\tilde{\pi})$, if the analytic conductors of $\pi$ and $\pi'$ are less than $Q$, then $Y(\pi,\pi')=O(\log^2 Q)$ 
(see, for example, [G-H]).

\qquad One wants to give an upper bound on $Y(\pi,\pi')$ which grows moderately in $Q$ without assuming a Riemann Hypothesis.  
In certain settings, this can be done through non-analytic means. As an example, Murty [Mu] used the Riemann-Roch theorem on the 
modular curve $X_0(N)$ to show that when $\pi$ and $\pi'$ correspond to holomorphic modular forms of level $N$ and even weight 
$k$, then $Y(\pi,\pi')=O(kN\log\log N)$.  For the case of Maass forms on the upper half plane, Huntley [H] used the method of 
Rayleigh quotients to show that $Y(\pi,\pi')$ grows at most linearly in the eigenvalue.  More recently, Baba, Chakraborty, and 
Petridis [B-C-P] proved a linear bound in the level and weight of holomorphic Hilbert modular forms, again using Rayleigh quotients.

\qquad This paper is concerned with, among other things, the determination of cusp forms on ${\rm GL}_n$ by their first few local 
components when measured with respect to both the archi-medean and non-archimedean parameters.  This case has been treated elsewhere 
by Moreno [Mo2], who derived a polynomial bound for $Y(\pi,\pi')$ when $n=2$ but could do no better than $Y(\pi,\pi')=O(e^{A\log^2 Q})$ 
for some constant $A>0$ when $n\geq 3$.  Moreno's idea was to demonstrate a region of non-vanishing for $L(s,\pi\times\pi')$ within 
the critical strip and apply this to an explicit formula relating sums over zeros to sums over primes.  For this strategy to work, 
quite a wide zero-free region is needed, one which decays logarithmically in all paramenters (with the possible exception of one real 
zero).  Unable to obtain this for $n$ greater than 2, Moreno used the phenomenon of zero repulsion to extract his exponential bound.  
In this paper, we obtain a modest zero-free region for $L(s,\pi\times\pi')$ for all $n\geq 2$, decaying {\it polynomially} in all 
parameters, and deduce from this, through an elementary method which, by contrast with Moreno's, uses sums over integers rather than 
primes, that $Y(\pi,\pi')=O(Q^A)$ for some constant $A>0$.  

\qquad Throughout this paper $\pi$ and $\pi'$ will denote (unitary) cuspidal automorphic representations of ${\rm GL}_n(\mathbb{A})$ 
($n\geq 1$).  We will make the implicit assumption that the central characters of $\pi$ and $\pi'$ are trivial on the product of positive 
reals $\mathbb{R}^+$ when embedded diagonally into the (archimedean places of) the ideles.  Under this normalization the Rankin-Selberg 
product $L(s,\pi\times\pi')$ has a pole at $s=1$ if and only if $\pi'=\tilde\pi$. 

\qquad The starting point of our inquiry is our Theorem \ref{polar part} where we give a lower bound on the polar part of 
$L(s,\Pi\times\widetilde\Pi)$ for $\Pi$ an isobaric representation of ${\rm GL}_d(\mathbb{A})$.  Theorem \ref{polar part} is proven 
through an approximation of the polar part by a smooth average of the coefficients of $L(s,\Pi\times\widetilde{\Pi})$.  That these 
coefficients are non-negative, and that certain of them are bounded away from zero (our Proposition \ref{Shay's save}), ensures that 
their average cannot be too small.  The error in the approximation is controlled through Mellin inversion by the functional equation 
of $L(s,\Pi\times\widetilde{\Pi})$ and is negligable as soon as the length of the sum is a large enough power of $Q$.

\qquad We then proceed to derive a first consequence of Theorem \ref{polar part}, proving an inverse polynomial lower bound on the 
edge of the critical strip for Rankin-Selberg $L$-functions $L(s,\pi\times\pi')$.  To simplify the statement, we write $Aut_n(\leq Q)$ 
for the set of all cuspidal automorphic representations $\pi$ of ${\rm GL}_n(\mathbb{A})$ with analytic conductor $\mathcal{C}(\pi)$ 
less than $Q$.

\bigskip

\qquad {\sc Theorem \ref{RS guy}.} {\it  Let $\pi\in Aut_n(\leq Q)$ and $\pi'\in Aut_{n'}(\leq Q)$ and assume that $\pi\neq\tilde\pi'$.  
Let $t\in\mathbb{R}$.  There exists $A=A(n,n')>0$ such that}

\begin{equation*}
|L(1+it,\pi\times\pi ' )|\underset{n,n'}{\gg} (Q(1+|t|))^{-A}.
\end{equation*}

\bigskip

\qquad To prove Theorem \ref{RS guy}, we apply Theorem \ref{polar part} to the isobaric sum 
$\Pi=\pi\otimes|\det|^{it/2}\boxplus\pi'\otimes|\det|^{it/2}$ on ${\rm GL}_d$, where $d=n+n'$.  With this choice of $\Pi$, we force 
the polar part of $L(s, \Pi\times\widetilde\Pi)$ to contain $L(1+it,\pi\times\pi')$ as a factor.  The convexity principle can be used 
to bound the factors that remain from above by a power of $Q$.  Since $\mathcal{C}(\Pi\times\widetilde\Pi)$ is itself bounded by a 
power of $Q$, we can then make the passage from the lower bound furnished by Theorem \ref{polar part} to that for $L(1+it,\pi\times\pi')$ 
as stated in Theorem \ref{RS guy}.

\qquad Theorem \ref{RS guy} has already found applications elsewhere in the literature.  For instance, Lapid [L] has recently 
shown that a lower bound on $L(1+it,\pi\times\pi')$ that decays at most polynomially in $Q(1+|t|)$ is a central issue in the 
convergence of Jacquet's relative trace formula.  

\qquad From Theorem \ref{RS guy} it is a short hop to obtain narrow zero-free regions.  It is known (see, for instance, [S]) that when both $\pi$ and $\pi'$ are self-dual, the method of de la Vall\'ee Poussin can be carried out successfully to give a (wide) zero-free region for $L(s,\pi\times\pi')$ of logarithmic 
type when the imaginary parameter $|t|\geq 1$.  When exactly one is self-dual, a standard zero-free region can be derived for all $t$.  Making the most of recent progress in functoriality, Ramakrishnan and Wang [R-W] have eliminated any
assumption of self-duality in certain low rank cases.   More precisely, they show that for $\pi$ and $\pi'$ on ${\rm GL}_2$ over 
$\mathbb{Q}$, the $L$-functions $L(s,\pi\times\pi')$ and $L(s, {\rm sym}^2\pi\times {\rm sym}^2\pi)$, as long as they are not divisible 
by $L$-functions of quadratic characters, admit no Seigel zeros.  For the cases that remain, we derive as a simple consequence of 
Theorem \ref{RS guy} a zero-free region for $L(s,\pi\times\pi')$ for arbitrary $\pi$ and $\pi'$, the width of which tapers polynomially 
in all parameters, and which remains valid even for $t=0$.  This is recorded in Corollary \ref{whodunnit?}.

\qquad The methods contained in Sections \ref{residue section} and \ref{lower bound section}, which combine to give Theorem \ref{RS guy}, 
can be thought of as an effectuation of Landau's lemma.  By contrast, Sarnak outlines a technique in [S] to show effective non-vanishing 
of $L$-functions through poles of Eisenstein series, and this too has now been carried out successfully by Gelbart, Lapid, and Sarnak 
[G-L-S].  These latter authors use the Langlands-Shahidi method to prove an inverse polynomial lower bound of certain $L$-functions 
along ${\rm Re}(s)=1$, but in the $t$-aspect only (and away from the real line).  Relative to the setting of our Theorem \ref{RS guy}, their result applies to a 
(at present) much larger class of $L$-functions.  Namely, to any $L$-function or product of $L$-functions obtained as the residue 
of an Eisenstein series they give a lower bound along ${\rm Re}(s)=1$; without full functoriality, it cannot be said that each one of 
these is the $L$-function of an automorphic form on ${\rm GL}_n$.  One striking application given by the authors of [G-L-S] is to 
$L(s,\pi, {\rm sym}^9)$, the symmetric-ninth power $L$-function of a cusp form $\pi$ on ${\rm GL}_2$: they prove a lower bound for 
$L(s,\pi,{\rm sym}^9)$ along ${\rm Re}(s)=1$ despite the fact that it is not yet known whether $L(s,\pi, {\rm sym}^9)$ is zero-free 
to the right of $1$. 

\qquad The final section in this paper is devoted to deriving the following effective multiplicity one statement.  In the proof, we 
exploit the fact that, with the aforementioned normalization on the central character, $L(s,\pi\times\tilde{\pi}')$ has a pole at 
$s=1$ if and only if $\pi=\pi'$.   The idea is that Theorem \ref{polar part} quantifies this property by providing a lower bound on 
the residue of $L(s,\pi\times\tilde{\pi})$.

\bigskip

\qquad {\sc Theorem \ref{mult one}.}  {\it Let $n\geq 1$.  Let $\pi=\otimes_v \pi_v$ and $\pi '=\otimes_v \pi_v '$ be in $Aut_n(\leq Q)$.  
Denote by $S$ the set of all finite places of $F$ at which either $\pi$ or $\pi'$ is ramified.  There exist constants $c=c(n)>0$ and 
$B=B(n)>0$ such that if $\pi_\textfrak{p}\simeq \pi_\textfrak{p}'$ for all prime ideals $\textfrak{p}\notin S$ with absolute norm 
${\rm N}\textfrak{p} \leq c Q^{B}$, then $\pi = \pi '$.}

\bigskip

\qquad The proof of Theorem \ref{mult one} allows for a weakening of the hypotheses, to the extent that one may suppose a mere 
approximate equivalence between the Dirichlet coefficients of the two forms and still retain the conclusion.  In this way we are 
able to deduce in Corollary \ref{finite} that the set $Aut_n(\leq Q)$ is finite. 

\medskip

\qquad {\it Acknowledgements:}  This paper came about through the suggestion of my thesis advisor, Peter Sarnak, and I would like to 
thank him now for his continual encouragement in this project.  I am also happy to acknowledge Erez Lapid for making important 
suggestions to correct the exposition and for explaining his own related work.  Lastly, I am indebted to Akshay Venkatesh for pointing 
out an improvement to an earlier version of Lemma \ref{Shay's save}.

\begin{center}
\section{\rm Preliminaries on $L$-functions}\label{L-functions}
\end{center}

In this section we give basic notation and definitions of standard and Rankin-Selberg $L$-functions, including their fundamental analytic properties and functional equations.

\bigskip

\quad {\it Standard $L$-function.}  Let $\pi$ be a cusp form on ${\rm GL}_n$ over a number field $F$.  To every prime ideal $\textfrak{p}$ at which $\pi_\textfrak{p}$ is unramified there is an associated set of $n$ non-zero complex Satake parameters $\{ \alpha_\pi (\textfrak{p},i)\}$ out of which one may define local $L$-functions

\begin{equation}\label{unramified standard}
L(s,\pi_\textfrak{p})=\prod_{i=1}^n (1-\alpha_\pi (\textfrak{p},i){\rm N}\textfrak{p}^{-s})^{-1}.
\end{equation}

At $\textfrak{p}$ where $\pi_\textfrak{p}$ is ramified the local $L$-function is defined in terms of the Langlands parameters of $\pi_\textfrak{p}$.  It is of the form $L(s,\pi_\textfrak{p})=P_\textfrak{p}({\rm N}\textfrak{p}^{-s})^{-1}$ where $P_\textfrak{p}(x)$ is a polynomial of degree at most $n$ and $P_\textfrak{p}(0)=1$.  It is possible to write the local factors at ramified primes in the form of (\ref{unramified standard}) with the convention that some of the $\alpha_\pi (\textfrak{p},i)$'s may be zero.  The $\alpha_\pi$ satisfy the bound

\begin{equation}\label{finite LRS bounds}
|\alpha_\pi (\textfrak{p},i)|\leq {\rm N}\textfrak{p}^{1/2-(n^2+1)^{-1}}
\end{equation}

by the work of Luo-Rudnick-Sarnak [L-R-S].

\qquad At each archimedean place $v$ a set of $n$ complex Langlands parameters $\{\mu_\pi (v,i)\}_{i=1}^n$ is associated to $\pi_v$.  The local factor at $v$ is defined to be

\begin{equation*}
L(s,\pi_v )=\prod_{i=1}^n \Gamma_{F_v}(s+\mu_\pi (v,i)),
\end{equation*}

where $\Gamma_{\mathbb{R}}(s)=\pi^{-s/2}\Gamma (s/2)$ and $\Gamma_{\mathbb C} (s)=2(2\pi)^{-s}\Gamma (s)$.  The $\mu_\pi$ satisfy 

\begin{equation}\label{infinite LRS bounds}
|{\rm Re} \ \mu_\pi (v,i)|\leq 1/2-(n^2+1)^{-1}
\end{equation}

again by [L-R-S].  

\qquad We denote by $\tilde{\pi}$ the contragredient representation of $\pi$. It is an irreducible cuspidal representation of ${\rm GL}_n(\mathbb{A})$.  For any place $v$ of $F$, $\tilde\pi_v$ is equivalent to the complex conjugate $\overline\pi_v$ [G-K], and hence

\begin{equation*}
\{ \alpha_{\tilde\pi}(\textfrak{p},i)\}=\{\overline{\alpha_\pi (\textfrak{p},i)}\}\qquad\text{and}\qquad \{\mu_{\tilde\pi}(v,i)\}=\{\overline{\mu_\pi (v,i)}\}.
\end{equation*}

\qquad By the bounds in (\ref{finite LRS bounds}), the product $\prod_{\textfrak{p}<\infty} L(s,\pi_\textfrak{p})$ converges absolutely on ${\rm Re}(s)>3/2-(n^2+1)^{-1}$ (in fact on ${\rm Re} (s)>1$, by Rankin-Selberg theory).  We write this product as a Dirichlet series over the integral ideals of the ring of integers $\mathcal{O}_F$ of $F$:

\begin{equation*}
L(s,\pi )=\prod_{\textfrak{p}<\infty } L(s,\pi_\textfrak{p} )=\sum_{\textfrak n} \lambda_\pi (\textfrak{n}){\rm N}\textfrak{n}^{-s}.
\end{equation*}

Let $S_\infty$ denote the set of the infinite places.  The complete $L$-function, defined to be $\Lambda (s, \pi )=L(s,\pi)\prod_{v\in S_\infty}L(s,\pi_v )$, is an entire function (except when $\pi$ is the trivial representation on ${\rm GL}_1$ so that $L(s,\pi)$ is the zeta function).  $\Lambda (s,\pi)$ has order 1 and is bounded in vertical strips.  It satisfies a functional equation $\Lambda (s,\pi )=W(\pi) q(\pi)^{1/2 - s}\Lambda (1-s,\tilde{\pi} )$ where $q(\pi)$ is the {\it arithmetic conductor} and $W (\pi)$, a complex number of modulus 1, is the {\it root number}.  We define

\begin{equation*}
\lambda_\infty (\pi; t)=\prod_{i=1}^n\prod_{v\in S_\infty}(1+|it+\mu_\pi (v, i)|)
\end{equation*}

and call $\mathcal{C}(\pi;t)=q(\pi)\lambda_\infty (\pi;t)$ the {\it analytic conductor} of $\pi$ (along the line $s=1+it$).  This definition was originally given in [I-S].  We denote $\mathcal{C}(\pi;0)$ by $\mathcal{C}(\pi)$.  

\bigskip

\quad {\it Rankin-Selberg $L$-functions.}  Let $\pi=\otimes_v \pi_v$ and $\pi '=\otimes_v \pi'_v$ be cuspidal representations of ${\rm GL}_n(\mathbb{A})$ and ${\rm GL}_{n'}(\mathbb{A})$.  For prime ideals $\textfrak{p}$ at which neither $\pi_\textfrak{p}$ nor $\pi_\textfrak{p}'$ is ramified let $\{\alpha_\pi (\textfrak{p},i)\}_{i=1}^n$ and $\{\alpha_{\pi'}(\textfrak{p},i)\}_{i=1}^{n'}$ be the respective Satake parameters of $\pi$ and $\pi'$.  The Rankin-Selberg $L$-function at such a $\textfrak{p}$ is defined to be

\begin{equation*}
L(s,\pi_\textfrak{p}\times\pi'_\textfrak{p})=\prod_{i=1}^n\prod_{j=1}^{n'} (1-\alpha_{\pi\times \pi'} (\textfrak{p},i,j){\rm N}\textfrak{p}^{-s})^{-1}.
\end{equation*}

These parameters satisfy

\begin{equation}\label{finite LRS RS bounds}
|\alpha_{\pi\times\pi'}(\textfrak{p},i,j)|\leq 1-(n^2+1)^{-1}-(n'^2+1)^{-1}.
\end{equation}

At primes at which either $\pi_\textfrak{p}$ or $\pi'_\textfrak{p}$ is unramified we have 

\begin{equation*}
\{\alpha_{\pi\times \pi'} (\textfrak{p},i,j)\}=\{\alpha_\pi (\textfrak{p},i)\alpha_{\pi'} (\textfrak{p},j)\}.
\end{equation*}

\qquad At each infinite place $v$ there exists a set of $nn'$ parameters $\{\mu_{\pi\times\pi'}(v,i,j)\ |\ 1\leq i\leq n, 1\leq j\leq n'\}$ such that the local factor at $v$ is

\begin{equation*}
L(s,\pi_v \times\pi'_v )=\prod_{i=1}^n\prod_{j=1}^{n'}\Gamma_{F_v}(s+\mu_{\pi\times\pi'}(v,i,j)).
\end{equation*}

At any place $v$ we have

\begin{equation}\label{conjugate RS}
\{\overline{\mu_{\pi\times\pi'}(v,i,j)}\}=\{\mu_{\tilde\pi\times\tilde\pi'}(v,i,j)\}
\end{equation}

and

\begin{equation}\label{infinite LRS RS bound}
|{\rm Re}\ \mu_{\pi\times\pi'}(v,i,j)|\leq 1-(n^2+1)^{-1}-(n'^2+1)^{-1}.
\end{equation}

When the infinite place $v$ is unramified for both $\pi$ and $\pi'$ we have 

\begin{equation}\label{unrammy}
\{\mu_{\pi\times\pi'}(v,i,j)\}=\{\mu_\pi(v,i)+\mu_{\pi'}(v,j)\}.
\end{equation}

\qquad By the bounds (\ref{finite LRS RS bounds}), the product $\prod_{\textfrak{p}<\infty} L(s,\pi_\textfrak{p}\times\pi'_\textfrak{p})$ converges absolutely in ${\rm Re}(s)> 2-(n^2+1)^{-1}-(n'^2+1)^{-1}$.  We write this product as a Dirichlet series over all integral ideals of the ring of integers $\mathcal{O}_F$ of $F$:

\begin{equation*}
L(s,\pi\times\pi' )=\prod_{\textfrak{p}<\infty}L(s,\pi_\textfrak{p}\times\pi'_\textfrak{p})= \sum_\textfrak{n} \lambda_{\pi\times\pi'} (\textfrak{n}){\rm N}\textfrak{n}^{-s}.
\end{equation*}

It can be shown through Rankin-Selberg integrals [J-PS-S] that the Euler product $L(s,\pi\times\pi')$ actually converges in ${\rm Re}(s)>1$.  With $S_\infty$ as usual representing the set of infinite places, the completed $L$-function $\Lambda (s,\pi\times\pi')=L(s,\pi\times\pi')\prod_{v\in S_\infty}L(s,\pi_v\times\pi'_v)$ extends to a meromorphic function on $\mathbb{C}$, is bounded (away from its poles) in vertical strips, and is of order 1.  Under our normalization on the central characters, $\Lambda (s,\pi\times\pi')$ is entire if and only if $\tilde\pi\neq\pi'$.  The poles of $\Lambda (s,\pi\times\tilde\pi)$ are simple and are located at $s=1$ and $s=0$.

\qquad The functional equation $\Lambda (s,\pi\times\pi ' )= W(\pi\times\pi ') q(\pi\times\pi')^{1/2-s}\Lambda (1-s, \tilde{\pi}\times\tilde{\pi}' )$ is valid for all $s$, where $q(\pi\times\pi')$ is the {\it arithmetic conductor} and $W (\pi\times\pi')$, a complex number of modulus 1, is the {\it root number}.  Let

\begin{equation*}
\lambda_\infty (\pi\times\pi';t)=\prod_{i=1}^n\prod_{j=1}^{n'}\prod_{v\in S_\infty}(1+|it+\mu_{\pi\times\pi'}(v,i,j)|).
\end{equation*}

As in [I-S] we define $\mathcal{C}(\pi\times\pi';t)=q(\pi\times\pi')\lambda_\infty(\pi\times\pi';t)$ to be the {\it analytic conductor} of the $L$-function $L(s,\pi\times\pi')$.  We write $\mathcal{C}(\pi\times\pi'):=\mathcal{C}(\pi\times\pi';0)$.

\bigskip

\quad {\it Separation of Components.}  We have $\lambda_\infty(\pi\times\pi';t)\leq \lambda_\infty (\pi;0)^{n'}\lambda_\infty (\pi';t)^n$.  For unramified places this is easy to see by (\ref{unrammy}).  For the ramified infinite places, see the calculations in [R-S, Appendix].  The arithmetic conductor $q(\pi\times\pi')$ separates according to the following result of Bushnell-Henniart [B-H]: $q(\pi\times\pi')\leq q(\pi)^{n'}q(\pi')^n/(q(\pi),q(\pi'))$.  These together produce

\begin{equation}\label{bush-henn}
\mathcal{C}(\pi\times\pi';t)\leq \mathcal{C}(\pi)^{n'}\mathcal{C}(\pi';t)^n\leq C(\pi)^{n'}C(\pi')^n(1+|t|)^{nn'[F:\mathbb{Q}]}.
\end{equation}

\bigskip

\quad {\it Preconvex bound.}  Let $\mu\in\mathbb{C}$ be such that ${\rm Re}\mu\geq -1+\theta$ for some $\theta>0$.  By Stirling's asymptotic formula for the Gamma function, for $s=\sigma+it$ where $\sigma<\theta$,

\begin{equation*}
\frac{\Gamma ((1-s+\overline{\mu})/2 )}{\Gamma ((s+\mu )/2)}\underset{\sigma}{\ll} \ (1+|it+\mu|)^{1/2-\sigma}.
\end{equation*}

Let $\theta=(n^2+1)^{-1}$ and $\theta'=(n'^2+1)^{-1}$.  When combined with the duplication formula $\Gamma_\mathbb{C}(s)=\Gamma_\mathbb{R}(s)\Gamma_\mathbb{R}(s+1)$ and displays (\ref{conjugate RS}) and (\ref{infinite LRS RS bound}), this gives the following estimate on the quotient for $\sigma< \theta+\theta'$:

\begin{equation}
\label{Gamma asymptotic bound}
\frac{L(1-s,\tilde{\pi}_v\times\tilde{\pi}'_v)}{L(s,\pi_v\times\pi'_v)}\underset{\sigma}{\ll}\lambda_\infty(\pi\times\pi';t)^{1/2-\sigma}.
\end{equation}  

\qquad From the bounds (\ref{finite LRS RS bounds}) we deduce $L(s,\pi\times\pi')=O(1)$ on ${\rm Re}(s)\geq\sigma_0$ for any $\sigma_0>2-\theta-\theta'$.  By the functional equation and the above estimate (\ref{Gamma asymptotic bound}), $L(\sigma+it,\pi\times\pi')=O(\mathcal{C}(\pi\times\pi';t)^{1/2-\sigma})$ on $\sigma\leq\sigma_0$ for any $\sigma_0< -1+\theta+\theta'$.  Using the Phragm\'en-Lindel\"of principle and the nice analytic properties of $L(s,\pi\times\pi')$ the following preconvex bound in the interval $-1+\theta+\theta'\leq \sigma\leq 2-\theta-\theta'$ is obtained:

\begin{equation}\label{first convexity}
L(\sigma+it,\pi\times\pi')\ll_\epsilon \mathcal{C}(\pi\times\pi';t)^{l(\sigma)+\epsilon},
\end{equation}

where $l(\sigma)$ is the linear function satisfying $l(-1+\theta+\theta')=3/2-\theta-\theta'$ and $l (2-\theta-\theta')=0$.  Note that the slope of $l(\sigma)$ is $-1/2$, regardless of $\theta, \theta'$.

\bigskip

\quad {\it Isobaric representations.}  An isobaric representation $\Pi$ on ${\rm GL}_d$ can be written 

\begin{equation*}
\Pi=\pi_1\otimes|\det|^{it_1}\boxplus\cdots\boxplus\pi_\ell\otimes|\det|^{it_\ell}, 
\end{equation*}

where $\pi_j$ is a cusp form on ${\rm GL}_{n_j}$ with $\sum_i n_i=d$, and $t_j\in\mathbb{R}$.  The $L$-function $L(s,\Pi)$ decomposes as a product $L(s,\Pi)=\prod_j L(s+it_j,\pi_j)$, and its analytic conductor is $\mathcal{C}(\Pi;t)=\prod_j\mathcal{C}(\pi_j;t+t_j)$.  Let 

\begin{equation*}
\Pi'=\pi'_1\otimes|\det|^{it'_1}\boxplus\cdots\boxplus\pi'_{\ell'}\otimes|\det|^{it'_{\ell'}}
\end{equation*}

be another isobaric representation on ${\rm GL}_{d'}$, with $\pi_j'$ on ${\rm GL}_{n_j'}$, $\sum_i n_i'=d'$, and $t'_j\in\mathbb{R}$.  Then the Rankin-Selberg product is $L(s,\Pi\times\Pi')=\prod_{j,k}L(s+t_j+t'_k,\pi_j\times\pi_k')$ with analytic conductor

\begin{equation}\label{factor con}
\mathcal{C}(\Pi\times\Pi';t)=\prod_{j=1}^\ell\prod_{k=1}^{\ell'}\mathcal{C}(\pi_j\times\pi_k';t+t_j+t'_k).
\end{equation}

As usual we set $\mathcal{C}(\Pi\times\Pi')=\mathcal{C}(\Pi\times\Pi';0)$.

\bigskip

\section{\rm A lower bound on the polar part of $L(s,\Pi\times\widetilde\Pi)$}\label{residue section}

\medskip

\qquad The goal of this section is to establish Theorem \ref{polar part} wherein a lower bound is given on the polar part of $L(s, \Pi\times\widetilde\Pi)$ for $\Pi$ an isobaric representation of ${\rm GL}_d(\mathbb{A})$.  We preface the proof by two lemmas.  Lemma \ref{Shay's save} shows that certain of the coefficients found in the Dirichlet series of $L(s, \Pi\times\widetilde\Pi)$ are bounded away from zero.  In Lemma \ref{second lemma}, this fact combines with the positivity of each one of the coefficients to bound their partial sum from below by a positive power of the length.  Theorem \ref{polar part} will then be shown to follow from these two lemmas through Mellin inversion.

\bl\label{Shay's save} Let $d\geq 1$.  For non-zero complex numbers $\alpha_1,\ldots \alpha_d$ define the coefficients $b_k$ by

\begin{equation}\label{d coeff}
\sum_{k\geq 0}b_k X^k=\prod_{i=1}^d \prod_{j=1}^d (1-\alpha_i\overline{\alpha_j}X)^{-1}.
\end{equation}

If the $\alpha$ satisfy $|\prod_{i=1}^d\alpha_i|=1$, then $b_d\geq 1$.\el

\pf A partition $\lambda=(\lambda_i)$ is a sequence of nonincreasing non-negative integers $\lambda_1\geq\lambda_2\geq\ldots$ with only finitely many non-zero entries.  For a partition $\lambda$, denote by $\ell(\lambda)$ the number of non-zero $\lambda_i$, and set $|\lambda|=\sum_i \lambda_i$.  For $\lambda$ such that $\ell (\lambda)\leq d$, let $s_\lambda (\alpha) $ be the Schur polynomial associated to $\lambda$, that is,

\begin{equation*}
s_\lambda (\alpha)={\rm det}(\alpha_i^{\lambda_j+d-j})_{ij}\ /\ {\rm det}(\alpha_i^{d-j})_{ij}.
\end{equation*}

By the orthogonality of the Schur polynomials (see, for instance, [Ma]),

\begin{equation*}
\prod_{i=1}^d \prod_{j=1}^d (1-\alpha_i\overline{\alpha_j}X)^{-1}=\sum_{\ell (\lambda)\leq d} |s_\lambda (\alpha)|^2X^{|\lambda|}.
\end{equation*}

For $\lambda=(\lambda_1,\ldots ,\lambda_d,0,\ldots )$, set $\hat\lambda=(\lambda_1-\lambda_d,\ldots ,\lambda_{d-1}-\lambda_d,0,\ldots )$.  Then $s_\lambda (\alpha)=\alpha_1^{\lambda_d}\cdots \alpha_d^{\lambda_d}s_{\hat\lambda}(\alpha)$, and since $|\prod_{i=1}^d\alpha_i|=1$, this gives $|s_\lambda (\alpha)|^2=|s_{\hat\lambda}(\alpha)|^2$.  Furthermore, for any pair $(\lambda, k)$, where $\lambda$ is a partition satisfying $\ell (\lambda)\leq d-1$ and $k\geq 0$ is an integer, there exists a unique partition $\lambda^{(k)}$ with $\ell (\lambda^{(k)})\leq d$ and $|\lambda^{(k)}|=|\lambda|+kd$ such that $\widehat{\lambda^{(k)}}=\lambda$.  This implies

\begin{equation*}
\sum_{\ell (\lambda)\leq d} |s_\lambda (\alpha)|^2X^{|\lambda|}=(1-X^d)^{-1}\sum_{\ell (\lambda)\leq d-1} |s_\lambda (\alpha)|^2 X^{|\lambda |}.
\end{equation*}

If $|\lambda |=0$ then $s_\lambda (\alpha)=1$.  The $d$-th coefficient in (\ref{d coeff}) is therefore

\begin{equation*}
b_d=1+\underset{\ell (\lambda)\leq d-1}{\underset{|\lambda|=d}{\sum}}|s_\lambda (\alpha)|^2.
\end{equation*}

From this we glean $b_d\geq 1$, as desired.\qed

\bigskip

\qquad Let $S$ be a finite set of prime ideals for the integer ring $\mathcal{O}_F$ of the number field $F$.  Write $\mathbb{S}=\prod_{\textfrak{p}\in S} \textfrak{p}$.  Let $d\geq 1$.  For each prime $\textfrak{p}\notin S$, let there be associated $d$ non-zero complex numbers $\alpha (\textfrak{p},1),\ldots ,\alpha (\textfrak{p}, d)$.  Let $b(\textfrak{n})$ be a sequence of non-negative real numbers indexed by the integral ideals of $\mathcal{O}_F$.  Assume that $b(\textfrak{1})=1$ and that for $\textfrak{p}\notin S$

\begin{equation}
\label{bi-variable D(s)}
\sum_{k\geq 0} b(\textfrak{p}^k)X^k=\prod_{i=1}^d \prod_{j=1}^d (1-\alpha(\textfrak{p},i)\overline{\alpha(\textfrak{p},j)}X)^{-1}.
\end{equation} 

\qquad Let $ \psi(x)\in C_c^\infty (0,\infty )$ be a non-negative function such that $\psi (x)\geq 1$ on $[1,2]$ and $\psi(0)=0$.  Let

\begin{equation*}
F(Y)=\sum_\textfrak{n}b(\textfrak{n})\psi ({\rm N}\textfrak{n}/Y),
\end{equation*}

the sum being over {\it all} integral ideals $\textfrak{n}$.  Since the coefficients $b(\textfrak{n})$ and $\psi$ itself are non-negative, it follows that $0\leq F(Y)$.  Had we chosen a smoothing function supported in an interval around 0, the identity $b(1)=1$ would further imply that $1\ll F(Y)$.  The following lemma enables us to to take $\psi(0)=0$, while still improving upon $1\ll F(Y)$ to show actual growth in the parameter $Y$ as soon as $Y$ is large enough.

\bl\label{second lemma} With the notation as above, there exists a constant $C=C(d)>0$ such that $F(Y)\gg Y^{1/d}(\log{Y})^{-1}$ for all $Y\gg_d (\log {\rm N}\mathbb{S})^C$.\el

\pf As the coefficients $b(\textfrak{n})$ and $\psi$ are non-negative, the sum $F(Y)$ can be truncated to give

\begin{equation*}
F(Y)\geq \sum_{Y\leq {\rm N}\textfrak{n}\leq 2Y}b(\textfrak{n})\psi ({\rm N}\textfrak{n}/Y)\geq \sum_{Y\leq {\rm N}\textfrak{n}\leq 2Y}b(\textfrak{n})\geq \sum_{Y\leq {\rm N}\textfrak{p}^d\leq 2Y} b(\textfrak{p}^d).
\end{equation*}

By (\ref{bi-variable D(s)}), the inequality $b(\textfrak{p}^d)\geq 1$ of Lemma \ref{Shay's save} may be applied to each $\textfrak{p}\notin S$.  Thus

\begin{align*}
F(Y) &\geq \#\{\textfrak{p}\ : \ Y^{1/d}\leq {\rm N}\textfrak{p}\leq (2Y)^{1/d}, \ \textfrak{p}\notin S\}\\
&\geq \#\{\textfrak{p}\ : \ Y^{1/d}\leq {\rm N}\textfrak{p}\leq (2Y)^{1/d}\} - \#\{\textfrak{p} : \textfrak{p}\in S\}:= A-B.
\end{align*}

As long as $A\geq 2B$ we have $F(Y)\geq \frac{1}{2}A$.  Since $B\leq\log {\rm N}\mathbb{S}$ and by the Prime Number Theorem $A\sim_d Y^{1/d}/\log Y$ (the implied constant depending also on the number field $F$), the lemma immediately follows.  \qed

\bigskip

\qquad Let $\ell$ be a positive integer and $\pi_i$, for $1\leq i\leq \ell$, be cuspidal automorphic representations of ${\rm GL}_{n_i}(\mathbb{A})$, $n_i\geq 1$.  For real numbers $t_1,\ldots ,t_\ell$ such that $t_i=t_j$ if $\pi_i=\pi_j$, let 

\begin{equation*}
\Pi=\pi_1\otimes|\det|^{it_1}\boxplus\cdots\boxplus\pi_\ell\otimes|\det|^{it_\ell}
\end{equation*}

be an isobaric representation on ${\rm GL}_d$, where $d=n_1+\cdots +n_\ell$.  The Rankin-Selberg $L$-function $L(s,\Pi\times\widetilde\Pi)$ has a pole of order at most $\ell^2$ at $s=1$ and, under our normalization on the central characters and the assumption on the twists $t_i$, is holomorphic elsewhere along ${\rm Re}(s)=1$.  Write the Laurent series expansion of $L(s,\Pi\times\widetilde\Pi)$ at $s=1$ as

\begin{equation*}
L(s,\Pi\times\widetilde\Pi)=\sum_{k=-\ell^2}^\infty r_k(s-1)^k.
\end{equation*}

\qquad The following theorem gives a lower bound on the polar part of $L(s,\Pi\times\widetilde\Pi)$ of polynomial decay in all parameters.  A result of this type was first proved by Carletti, Monte Bragadin and Perelli ([C-MB-P], Theorem 5).  Their approach, expressed in the language of Selberg class $L$-functions, uses both the positivity of the coefficients $b(\textfrak{n})$ and the identity $b(\textfrak{1})=1$.  These two data are therefore enough to buy a polynomial dependence on the conductor.  By incorporating the extra information contained in our Lemma \ref{Shay's save}, however, we improve the power of the conductor given by their technique.  Indeed without Lemma \ref{Shay's save} the lower bound in Theorem \ref{polar part} would be $\mathcal{C}(\Pi\times\Pi)^{-\frac{1}{2}+\epsilon}$.  It should be noted that Theorem \ref{polar part} in fact interpolates the bound $q^{-1/2-\epsilon}\ll_\epsilon L(1,\chi)$ for real primitive Dirichlet characters, making it a close approximation of Dirichlet's bound.  We discuss this in more detail in Example \ref{example} following the proof.

\bt\label{polar part} With the notation as above, for every $\epsilon>0$

\begin{equation*}
 \sum_{k=1}^{\ell^2} |r_{-k}|\underset{\epsilon}{\gg} \mathcal{C}(\Pi\times\Pi)^{-\frac{1}{2}(1-1/d)-\epsilon}.
\end{equation*}

\et

\pf Let $\psi(x)$ be a smooth compactly supported non-negative function on the positive reals with $\psi(x)\geq 1$ on $[1,2]$ and $\psi(0)=0$.  The Mellin transform of $\psi$,

\begin{equation*}
\hat{\psi}(s)=\int_0^\infty \psi(x)x^{s-1}dx
\end{equation*}

is an entire function with rapid decay in vertical strips.  As a Dirichlet series $L(s,\Pi\times\widetilde\Pi)$ can be written $L(s,\Pi\times\widetilde\Pi)=\sum_\textfrak{n} b(\textfrak{n}){\rm N}\textfrak{n}^{-s}$.  Let $F(Y)=\sum_\textfrak{n} b(\textfrak{n})\psi ({\rm N}\textfrak{n}/Y)$.  From the Mellin inversion formula it follows that

\begin{equation*}
F(Y)=\sum_{\textfrak{n}}b(\textfrak{n})\left( \frac{1}{2\pi i}\underset{\sigma =2}{\int} \hat{\psi}(s)(Y/{\rm N}\textfrak{n})^s\ ds  \right).
\end{equation*}

The absolute convergence of $L(s,\Pi\times\widetilde\Pi)$ beyond $\sigma=1$ allows us to switch the order of the sum and integral to obtain
 
\begin{equation}
\label{sum as integral}
F(Y)= \frac{1}{2\pi i}\underset{\sigma=2}{\int} L(s,\Pi\times\widetilde\Pi)\hat{\psi}(s)Y^s\ ds.
\end{equation}

The integrand in (\ref{sum as integral}) is bounded in vertical strips.  The principle of Phragm\'en-Lindel\"of thus allows the contour of integration to be shifted to the left, while picking up the residue of the integrand at $s=1$.  Shifting to the line $\sigma=-b$ for $b\geq 1$ we get

\begin{equation*}
F(Y)=\underset{s=1}{{\rm Res}}\ \hat\psi (s)L(s,\Pi\times\widetilde\Pi)Y^s+\frac{1}{2\pi i}\underset{\sigma=-b }{\int}L(s,\Pi\times\widetilde\Pi)\hat{\psi} (s)Y^s\ ds.
\end{equation*}

To estimate the integral, we use (\ref{first convexity}), (\ref{factor con}) and (\ref{bush-henn}), noting the rapid decay in vertical strips of the integrand, to obtain

\begin{equation}\label{integral upper bound}
F(Y)=\underset{s=1}{\rm Res}\ \hat\psi (s)L(s,\Pi\times\widetilde\Pi)Y^s+O_\epsilon (\mathcal{C}(\Pi\times\Pi)^{l(-b)+\epsilon}Y^{-b}).
\end{equation}

When $S$ is the set of primes at which $\Pi$ is ramified, the sum $F(Y)$ satisifies the conditions of Lemma \ref{second lemma}, so that for $Y\gg_d (\log\mathcal{C}(\Pi\times\Pi))^C$

\begin{equation}\label{compare}
Y^{1/d}(\log Y)^{-1}\ll_d F(Y).
\end{equation}

If we take $Y= c\ \mathcal{C}(\Pi\times\Pi)^{\frac{l(-b)+1}{b}}$ for a large enough constant $c>0$ then (\ref{compare}) is valid and the lower bound on $F(Y)$ in (\ref{compare}) dominates the error term in (\ref{integral upper bound}).  Given any $\epsilon>0$ we may take $b$ large enough with respect to $d$, $\epsilon$, and the constant term in $l(b)$ to ensure that $Y=\mathcal{C}(\Pi\times\Pi)^{1/2+\epsilon}$ is a stronger condition than that above (recall that the slope of $l(b)$ is $-1/2$).  With this value of $Y$ we obtain

\begin{equation}\label{residue0}
Y^{1/d}(\log Y)^{-1}\ll_d \underset{s=1}{\rm Res}\ \hat\psi (s)L(s,\Pi\times\widetilde\Pi)Y^s.
\end{equation}

Since $Y^s=Y\sum_{j\geq 0}(\log Y)^j(s-1)^j/j!$, the right hand side of (\ref{residue0}) is

\begin{equation}
\underset{\psi}{\ll} Y \sum_{j+k=-1} |r_{k}|(\log{Y})^j/j! \ll Y(\log Y)^{\ell^2-1} \sum_{k=1}^{\ell^2} |r_{-k}|.\label{residue1}
\end{equation}

The theorem follows upon substituting $Y=\mathcal{C}(\Pi\times\Pi)^{1/2+\epsilon}$ into (\ref{residue0}) and (\ref{residue1}).\qed

\bex\label{example} {\rm Let $\pi_1=\chi$, a primitive real Dirichlet character of modulus $q$, and $\pi_2=1$, the trivial character.  Then for $\Pi=\pi_1\boxplus\pi_2=\chi\boxplus 1$, we have}

\begin{equation*}
L(s,\Pi\times\widetilde\Pi)=[\zeta (s)L(s,\chi)]^2,
\end{equation*}

{\rm $\mathcal{C}(\Pi\times\Pi)\asymp q^2$ and $d=2$.  The function $L(s,\Pi\times\tilde\Pi)$ has a double pole at $s=1$ and nowhere else, and if we denote by $\gamma=\zeta'(1)$ Euler's constant, then}

\begin{equation*}
r_{-2}=L(1,\chi)^2 \quad\text{\rm and}\quad r_{-1}=2L'(1,\chi)L(1,\chi) + 2\gamma L(1,\chi)^2.
\end{equation*}

{\rm Applying Theorem \ref{polar part} gives}

\begin{equation*}
\frac{1}{q^{1/2+\epsilon}}\ll_\epsilon L(1,\chi)\bigg(L(1,\chi)(1+2\gamma)+2|L'(1,\chi)|\bigg).
\end{equation*}

{\rm Since $L^{(k)}(1,\chi)\ll_\epsilon (\log q)^k$, we conclude by this technique that}

\begin{equation*}
\frac{1}{q^{1/2+\epsilon}}\ll_\epsilon L(1,\chi),
\end{equation*}

{\rm which is only slightly worse than what Dirichlet deduced by his class number formula, namely $q^{-1/2}\ll L(1,\chi)$.}\eex

\begin{center}
\section{\rm Lower bounds for $L(1+it,\pi\times\pi')$}\label{lower bound section}
\end{center}

\qquad We shall now use Theorem \ref{polar part} to bound from below the value along ${\rm Re}(s)=1$ of the Rankin-Selberg $L$-function $L(s,\pi\times\pi')$.  To do so at the point $s=1+it$, we construct an auxillary Dirichlet series $L(s,\Pi \times\widetilde{\Pi})$ whose polar part contains $L(1+it,\pi\times\pi')$ as a factor.  Roughly speaking, this coincidence is ensured as soon as the order of the pole at $s=1$ is equal to the power to which $L(s,\pi\times\pi')$ divides $L(s,\Pi\times\widetilde\Pi)$.  This is precisely the case in which one classically appeals to Landau's lemma to show mere non-vanishing on the line ${\rm Re}(s)=1$.  The following theorem, Theorem \ref{RS guy}, can therefore be interpreted as an effectuation of Landau's lemma.

\qquad Now that in Example \ref{example} we have measured the quality of the exponent given by Theorem \ref{polar part}, we shall no longer give specific powers of the conductor in our results.  One reason for doing so is that the statements that follow all employ the preconvex bound (\ref{first convexity}) which can be improved by progress toward the Ramanujan conjecture (see [Molt]).  Beyond that subconvex bounds would improve the exponents even further.  There is therefore no compelling reason to specify each exponent, and we greatly simplify the exposition by not doing so.

\bigskip

\qquad {\sc Definition.}  For a real parameter $Q\geq 2$ we denote by  $Aut_n(\leq Q)$ the set of all cuspidal representations $\pi$ of ${\rm GL}_n(\mathbb{A})$ with analytic conductor $\mathcal{C}(\pi)$ less than $Q$.

\bigskip

\bt\label{RS guy} Let $\pi\in Aut_n(\leq Q)$ and $\pi'\in Aut_{n'}(\leq Q)$, and assume $\pi'\neq\tilde\pi$.  Let $t\in\mathbb{R}$.  There exists $A=A(n, n')>0$ such that

\begin{equation*}
|L(1+it,\pi\times\pi')|\underset{n,n'}{\gg} (Q(1+|t|))^{-A}.
\end{equation*}

\et

\pf Consider the unitary isobaric sum $\Pi=\pi\otimes|\det|^{it/2}\boxplus\ \pi'\otimes |\det|^{it/2}$, defined on ${\rm GL}_d$ where $d=n+n'$.  The Rankin-Selberg product $L(s,\Pi\times\widetilde{\Pi })$ can be written

\begin{equation*}
L(s,\pi\times\tilde{\pi})L(s,\pi '\times\tilde{\pi}')L(s+it,\pi\times\tilde{\pi} ')L(s-it,\tilde{\pi}\times\pi ' ).
\end{equation*}

We apply Theorem \ref{polar part} with $d=n+n'$ to get

\begin{equation}
\label{downers II'}
|r_{-1}|+|r_{-2}|\underset{\epsilon}{\gg} \mathcal{C}(\Pi\times\Pi)^{-\frac{1}{2}(1-1/d)-\epsilon}.
\end{equation}

By the factorization (\ref{factor con}) and the separation of components in (\ref{bush-henn}), the analytic conductor $\mathcal{C}(\Pi\times\Pi)$ of $L(s,\Pi \times\widetilde{\Pi})$ is 

\begin{equation*}
\mathcal{C}(\pi\times\pi)\mathcal{C}(\pi'\times\pi')\mathcal{C}(\pi\times\pi';t)^2\leq (1+|t|)^{2nn'[F:\mathbb{Q}]}Q^{4(n+n')}.
\end{equation*}

Thus the lower bound in (\ref{downers II'}) becomes

\begin{equation}\label{downers II}
|r_{-1}|+|r_{-2}|\gg (Q(1+|t|))^{-A_1}
\end{equation}

for some explicitly given $A_1=A_1(n,n')>0$.

\qquad Using $\overline{L(s,\pi\times\tilde\pi')}=L(\overline{s},\tilde\pi\times\pi')$ we compute $r_{-2}=R_{-1}R_{-1}'|L(1+it,\pi\times\tilde\pi')|^2$ and

\begin{align*}
r_{-1}=( R_{-1}R_0'+R_0R_{-1}') &|L(1+it,\pi\times\tilde\pi')|^2\\
&+2R_{-1}R_{-1}'{\rm Re} (L'(1+it,\pi\times\tilde\pi')\overline{L(1+it,\pi\times\tilde \pi')}).
\end{align*}

The inequality ${\rm Re}(z_1\overline{z_2})\leq |z_1z_2|$ and the preconvex bounds

\begin{equation*}
R_{-1}, R_{0}\ll Q^{A_2},\qquad L^{(k)}(1+it,\pi\times\tilde\pi')\ll (Q(1+|t|))^{A_2},
\end{equation*}

for $k=0, 1$, and some $A_2=A_2(n,n')>0$, now give

\begin{equation*}
r_{-1}, r_{-2}\ll |L(1+it,\pi\times\tilde\pi')|(Q(1+|t|))^{3A_2}.
\end{equation*}

When combined with (\ref{downers II}) this implies the theorem, the power being $A=A_1+3A_2$.\qed

\bigskip

\qquad As was mentioned in the introduction, the method of de la Vall\'ee Poussin can be used under certain cirmustances to derive zero-free regions for $L(s,\pi\times\pi')$ of logarithmic type.  For instance, to eliminate the possibility of a real exceptional zero for $L(s,\pi\times\pi')$, exactly one of $\pi$ and $\pi'$ must be self-dual.  In certain cases of low rank, Ramakrishnan and Wang [R-W] have eliminated the hypothesis of self-duality.  They show that for $\pi$ and $\pi'$ on ${\rm GL}_2$ over $\mathbb{Q}$, the $L$-functions $L(s,\pi\times\pi')$ and $L(s, {\rm sym}^2\pi\times {\rm sym}^2\pi)$, as long as they are not divisible by $L$-functions of quadratic characters, admit no Seigel zeros.  In all the cases that remain, the following corollary to Theorem \ref{RS guy} provides a healthy compromise. 

\bc\label{whodunnit?} Let $\pi\in Aut_n(\leq Q)$, $Aut_{n'}(\leq Q)$, and $t\in\mathbb{R}$.  There exist constants $c=c(n,n')>0$ and $A'=A'(n,n')>0$ such that $L(\sigma+it,\pi\times\pi')$ has no zeros in the interval

\begin{equation*}
1-\frac{c}{(Q(1+|t|))^{A'}}\leq\sigma\leq 1.
\end{equation*}

\ec

\bigskip

\pf Let $\beta+it$ denote the first zero of $L(s,\pi\times\pi')$ to the left of $1$ along the segment $\sigma+it$, $1/2<\sigma<1$.  Then we have

\begin{equation*}
L(1+it,\pi\times\pi')=\int_\beta^1 L'(\sigma+it,\pi\times\pi') \ d\sigma=(1-\beta)L'(\sigma_0+it,\pi\times\pi'),
\end{equation*}

for some $\beta\leq \sigma_0\leq 1$, by the mean value theorem.  We apply the preconvex bound for $L'(s,\pi\times\pi')$ on the critical line $\sigma=1/2$

\begin{equation*}
|L'(\sigma_0+it,\pi\times\pi')|\leq |L'(1/2+it,\pi\times\pi')|\ll (Q(1+|t|))^{A_3},
\end{equation*}

for some $A_3=A_3(n,n')>0$.  We finally apply the lower bound for $L(1+it,\pi\times\pi')$ from Theorem \ref{RS guy} to obtain the corollary, the power being $A'=A+A_3$.  \qed

\begin{center}
\section{\rm Effective multiplicity one}
\end{center}

\qquad We note that by Theorem \ref{polar part}, when $\pi\in Aut_n (\leq Q)$, we have

\begin{equation}\label{juicy!}
R:=\underset{s=1}{\rm Res}\ L(s,\pi\times\tilde\pi)\gg Q^{-B_1}
\end{equation}

for a constant $B_1=B_1(n)>0$.  

\bigskip

\bt\label{mult one} For a real parameter $Q\geq 1$ let $\pi,\pi'$ be in $Aut_n(\leq Q)$ and $S$ be any finite set of finite places of $F$ satisfying $|S|\ll \log Q$.  There exists a constant $B=B(n, S)>0$ such that if $\pi_\textfrak{p}\simeq\pi_\textfrak{p}'$ for all primes ideals $\textfrak{p}\notin S$ with ${\rm N} \textfrak{p} \leq Q^B$, then $\pi =\pi '$.

\et

\pf Fix as a test function any non-negative $\psi (x)\in C_c^\infty (0,\infty)$ with $\hat{\psi}(1)=1$.  Put $\mathbb{S}=\prod_{\textfrak{p}\in S}\textfrak{p}$ and define

\begin{equation*}
F(Y;\pi\times\tilde\pi')=\sum_{(\textfrak{n},\mathbb{S})=1} \lambda_{\pi\times\tilde\pi'}(\textfrak{n})\psi ({\rm N}\textfrak{n}/Y).
\end{equation*}

The hypothesis on the local representations means that the Satake parameters $\{\alpha_\pi(\textfrak{p},i)\}$ and $\{\alpha_{\pi'} (\textfrak{p},i)\}$ agree (as sets) for all prime ideals $\textfrak{p}\notin S$ with absolute norms within the specified range.  It follows that $\lambda_{\pi\times\tilde\pi}(\textfrak{p}^k)=\lambda_{\pi\times\tilde \pi'}(\textfrak{p}^k)$ for all primes ideals $\textfrak{p}\notin S$ with ${\rm N}\textfrak{p}\leq Q^B$ and all $k\geq 1$.  By multiplicativity on coprime ideals, one derives the condition that $\lambda_{\pi\times\tilde\pi}(\textfrak{n})=\lambda_{\pi\times\tilde \pi'}(\textfrak{n})$ for all ideals ${\rm N}\textfrak{n}\leq Q^B$ with $(\textfrak{n},\mathbb{S})=1$ -- that is 

\begin{equation}\label{new condition}
F(Y;\pi\times\tilde\pi)=F(Y;\pi\times\tilde \pi')\qquad\text{for}\qquad Y\leq Q^B.
\end{equation}

This will henceforth be our assumption.

\qquad With $S$ as in the statement of the theorem, let $L_S(s,\pi\times\tilde\pi')=\prod_{\textfrak{p}\notin S}L(s,\pi_\textfrak{p}\times\tilde\pi_\textfrak{p}')$ and $L^S(s,\pi\times\tilde\pi')=\prod_{\textfrak{p}\in S}L(s,\pi_\textfrak{p}\times\tilde\pi_\textfrak{p}')$.  Mellin inversion gives

\begin{equation*}
F(Y;\pi\times\tilde \pi')=\frac{1}{2\pi i}\underset{\sigma= 2}{\int}L_S(s,\pi\times\tilde\pi')\hat{\psi} (s)Y^s\ ds.
\end{equation*}

Let $\theta=\theta(n)=(n^2+1)^{-1}$, the quantity appearing in the Luo-Rudnick-Sarnak bounds (\ref{finite LRS bounds}).  We note that the local factor $L(s,\pi_\textfrak{p}\times\tilde\pi'_\textfrak{p})$, and thus the product $L^S(s,\pi\times\tilde\pi')$, is well-defined and invertible on ${\rm Re}(s)>1-2\theta$.  Since $L_S(s,\pi\times\tilde\pi')=L(s,\pi\times\tilde\pi')L^S(s,\pi\times\tilde\pi')^{-1}$, the first factor extending meromorphically to $\mathbb{C}$, we may move the contour to the line ${\rm Re}(s)=1-\theta$, while picking up the residue of the integrand at $s=1$.  This gives

\begin{equation*}
F(Y;\pi\times\tilde \pi')=\delta_{\pi,\pi'}YRL^S(1,\pi\times\tilde\pi)^{-1}+\frac{1}{2\pi i}\underset{\sigma= 1-\theta}{\int}L(s,\pi\times\tilde\pi')L^S(s,\pi\times\tilde\pi')^{-1}\hat{\psi} (s)Y^s\ ds.
\end{equation*}

\qquad We bound the individual factors in the above integrand.  The preconvex bound on $L(s,\pi\times\tilde \pi')$ at ${\rm Re}(s)=1-\theta$ is $L(1-\theta+it,\pi\times\tilde \pi')\ll (Q(1+|t|))^{B_2}$ for some $B_2=B_2(n, \theta)>0$.  By (\ref{finite LRS bounds}) we have for ${\rm Re}(s)=1-\theta$ 

\begin{equation*}
\prod_{1\leq i,j\leq n}|1-\alpha_\pi (\textfrak{p},i)\overline{\alpha_{\pi'}(\textfrak{p},j)}p^{-s}|\leq (1+p^{-\theta})^{n^2}=O_n(1).  
\end{equation*}

Since $|S|\ll\log Q$, this gives $|L^S(s,\pi\times\tilde\pi')|^{-1}=O(1)^{|S|}\leq Q^{B_3}$ for some constant $B_3=B_3(n, \theta)>0$.  By the rapid decay of $\hat{\psi}(s)$ along vertical lines then

\begin{equation}\label{integral uppers}
F(Y;\pi\times\tilde \pi')=\delta_{\pi,\pi'}YRL^S(1,\pi\times\tilde\pi)^{-1}+O (Y^{1-\theta}Q^{B_2+B_3}).
\end{equation}

\qquad Let $B>0$ be a constant such that (\ref{new condition}) holds and suppose that $\pi\neq\pi'$.  We seek a contradiction to the latter supposition.  The key observation is that under both (\ref{new condition}) and $\pi\neq\pi'$ the error term of $F(Y;\pi\times\tilde\pi)$ in equation (\ref{integral uppers}) must dominate the main term.  In this range, therefore, 

\begin{equation}\label{kitty}
Y^{\theta}=O(R^{-1}L^S(1,\pi\times\tilde\pi)Q^{B_2+B_3}).  
\end{equation}

Since $L(s,\pi\times\tilde\pi)$ has positive coefficients as a Dirichlet series to the right of 1, we can bound $L^S(1,\pi\times\tilde\pi)$ by the preconvex bound at $s=1$ of the regularization of $L(s,\pi\times\tilde\pi)$, so that $L^S(1,\pi\times\tilde\pi)=O(Q^{B_4})$.  By $R^{-1}\ll Q^{B_1}$ of display (\ref{juicy!}), equation (\ref{kitty}) becomes $Y=O( Q^{\theta^{-1}(B_1+B_2+B_3+B_4)})$.  To force a contradiction, we have only to take $B$ to be $B>\theta^{-1}(B_1+B_2+B_3+B_4).$ \qed

\brmk\label{relax} {\rm As we have seen, the condition of Theorem \ref{mult one} that the first few local components be isomorphic can be expressed instead as an equality of the initial coefficients of the Rankin-Selberg $L$-series.  In fact this latter condition can be relaxed to an approximate equivalence, in which the difference between the first few coefficients is bounded below by some expression in the conductor.  

\qquad Having chosen $\pi, \pi'\in Aut_n (\leq Q)$, let the set $S$ consist of precisely those prime ideals at which either $\pi$ or $\pi'$ is ramified.  Then $|S|\ll \log Q$ as required in the statement of Theorem \ref{mult one}.  Put $\mathbb{S}=\prod_{\textfrak{p}\in S}\textfrak{p}$.  We claim that if $\pi\neq\pi'$ then there exist numbers $B, C>0$ such that $|\lambda_\pi(\textfrak{n}_0)-\lambda_{\pi'}(\textfrak{n}_0)|\gg Q^{-C}$ for some square-free ideal $(\textfrak{n}_0, \mathbb{S})=1$ with ${\rm N}\textfrak{n}_0\leq Q^B$.  This relaxation is essential for comparing automorphic forms whose coefficients are not algebraic, as is believed to be the case for Maass wave forms.

\qquad By the previous arguments, since $\pi\neq\pi'$,}

\begin{equation*}
\sum_{(\textfrak{n},\mathbb{S})=1}(\lambda_{\pi\times\tilde\pi}(\textfrak{n})-\lambda_{\pi\times\tilde \pi'}(\textfrak{n}))\psi ({\rm N}\textfrak{n}/Y)=YRL^S(1,\pi\times\tilde\pi)^{-1}+O(Y^{1-\theta}Q^{B_2+B_3}).
\end{equation*}

{\rm Under Theorem \ref{mult one} if $Y=Q^B$ for $B$ large enough then it is the main term that dominates the error term, giving}

\begin{equation*}
\frac{1}{Y}\sum_{(\textfrak{n},\mathbb{S})=1}|\lambda_{\pi\times\tilde\pi}(\textfrak{n})-\lambda_{\pi\times\tilde \pi'}(\textfrak{n})|\psi({\rm N}\textfrak{n}/Y)\gg RL^S(1,\pi\times\tilde\pi)^{-1}\gg Q^{-B_1-B_4}.
\end{equation*}

{\rm There therefore exists an integral ideal $\textfrak{n}_0$ relatively prime to $\mathbb{S}$ with ${\rm N}\textfrak{n}_0\leq Q^B$ such that $|\lambda_{\pi\times\tilde\pi}(\textfrak{n}_0)-\lambda_{\pi\times\tilde \pi'}(\textfrak{n}_0)|\gg Q^{-B_1-B_4}$ (the implied constant depending on $\psi$).  At the cost of increasing $B_1, \ldots , B_4$ (and hence $B$ as well), the ideal $\textfrak{n}_0$ can be taken to be square-free (simply redo the proof of Theorem \ref{mult one} using the square-free unramified $L$-function).  Recall that $\lambda_{\pi\times\pi'}(\textfrak{n})=\lambda_\pi(\textfrak{n})\lambda_{\pi'}(\textfrak{n})$ on square-free unramified ideals $\textfrak{n}$.  By the bounds (\ref{finite LRS bounds}) with $\theta=(n^2+1)^{-1}$}

\begin{equation*}
|\lambda_{\pi\times\tilde\pi}(\textfrak{n}_0)-\lambda_{\pi\times\tilde \pi'}(\textfrak{n}_0)|\ll Q^{(1/2-\theta)B}|\lambda_\pi (\textfrak{n}_0)-\lambda_{\pi'}(\textfrak{n}_0)|
\end{equation*}

{\rm and the claim follows with $C=B_1+B_4+(1/2-\theta)B$. }\ermk

\bc\label{finite} The set $Aut_n(\leq Q)$ is finite.\ec

\pf Put $S=\{ \textfrak{p}\ : \ {\rm N}\textfrak{p}\leq Q\}$ and observe that the prime ideals at which any $\pi\in Aut_n(\leq Q)$ is ramified are contained in $S$.  Let $S=\bigsqcup S_i$ be a disjoint covering of $S$ by subsets $S_i$ satisfying $\prod_{\textfrak{p}\in S_i}{\rm N}\textfrak{p}\leq Q$.  Denote by $Aut_n(S_i)$ the set of all automorphic forms on ${\rm GL}_n/F$ unramified at finite places outside of $S_i$.  We have $Aut_n(\leq Q)\subset\bigsqcup Aut_n(S_i)$.  We shall show that each intersection $Aut_n (S_i)\cap Aut_n(\leq Q)$ is finite.

\qquad Let $B>0$ be a constant (to be fixed later).  Put $\mathbb{S}_i=\prod_{\textfrak{p}\in S_i}\textfrak{p}$.  For each $i$ let $\mathcal{I}_i$ the set of square-free ideals $(\textfrak{n}, \mathbb{S}_i)=1$ with ${\rm N}\textfrak{n}\leq Q^B$.  For constants $\epsilon, c>0$, consider the space of sequences of complex numbers

\begin{equation*}
X_i=X_i(\epsilon, c)=\{(\lambda (\textfrak{n}))_{\textfrak{n}\in\mathcal{I}_i}\ :\ |\lambda (\textfrak{n})|\leq c{\rm N}\textfrak{n}^{1/2-(n^2+1)^{-1}+\epsilon}\} 
\end{equation*}

endowed with the natural topology and metric as a closed subset of $\mathbb{C}^{M_i}$, where $M_i=|\mathcal{I}_i|$.  By the bounds (\ref{finite LRS bounds}), for any $\epsilon>0$ there exists a constant $c=c(\epsilon)>0$ such that the set $Aut_n(S_i)$ maps to $X_i$ via the Fourier coefficient map $FC_i: \pi\mapsto (\lambda_\pi (\textfrak{n}))_{\textfrak{n}\in\mathcal{I}_i}$.  Since $|S_i|\leq \log Q$ we may take $B$ as in Theorem \ref{mult one} to conclude that the restriction of $FC_i$ to $Aut_n (S_i)\cap Aut_n(\leq Q)$ is injective.  Moreover, the distance squared between any two $\pi$, $\pi'\in Aut_n(S_i)\cap Aut_n(\leq Q)$, considered as points in $X_i$, is

\begin{equation*}
dist (\pi,\pi')^2:=|FC(\pi)-FC(\pi')|^2=\sum_{\textfrak{n}\in\mathcal{I}_i}|\lambda_\pi (\textfrak{n})-\lambda_{\pi'}(\textfrak{n})|^2.
\end{equation*}

For $\pi$ and $\pi'$ distinct we thus have 

\begin{equation}\label{last one}
dist (\pi, \pi')\geq \underset{\textfrak{n}\in\mathcal{I}_i}{\rm max}\ |\lambda_\pi (\textfrak{n})-\lambda_{\pi'}(\textfrak{n})|\gg Q^{-C} 
\end{equation}

by Remark \ref{relax}.  Hence $Aut_n(S_i)\cap Aut_n(\leq Q)$ is discrete in $X_i$.  As $X_i$ is compact, the result follows. \qed

\brmk {\rm The bound on $|Aut_n(\leq Q)|$ given by the above corollary is probably very poor, possibly exponential.  Even though only $O(Q/\log Q)$ sets $S_i$ are needed to cover $S$, it is not evident that the lower bound (\ref{last one}) should be sufficient to prove that the slice $|Aut_n(S_i)\cap Aut_n(\leq Q)|$ itself is polynomial in $Q$.  A more sophisticated analysis using the trace formula should however give a sharp polynomial bound in all parameters.}\ermk

\begin{center}
{\sc References}
\end{center}

[B-H] C.J. Bushnell and G. Henniart. An upper bound on conductors for pairs. {\it J. Number Theory}  {\bf 65}  (1997),  no. 2, 183--196.

\bigskip

[B-C-P] S. Baba, K. Chakraborty, and Y. Petridis. On the number of Fourier coefficients that determine a Hilbert modular form.  {\it Proc. Amer. Math. Soc.} {\bf 130} (2002), no.9, 2497-2502.

\bigskip

[C-MB-P] E. Carletti, G. Monti Bragadin, and A. Perelli. On general $L$-functions.  {\it Acta Arith.}  {\bf 66}  (1994),  no. 2, 147--179.

\bigskip

[G-H] D. Goldfeld and J. Hoffstein. On the number of Fourier coefficients that determine a modular form.  {\it A tribute to Emil Grosswald: number theory and related analysis,}  385--393, Contemp. Math., 143, {\it Amer. Math. Soc., Providence, RI,} 1993.

\bigskip

[G-K] I.M. Gelfand and D. Kazhdan, ``Representation of the group ${\rm GL}(n,K)$ where $K$ is a local field'' in {\it Lie Groups and their Representations}, ed. by I.M. Gelfand, Wiley, New York, 1974, 95-118.

\bigskip

[G-L-S] S. Gelbart, E. Lapid, and P. Sarnak.  A new method for lower bounds of L-functions, {\it C. R. Acad. Sci.} Paris I 339 (2004) 91-94.

\bigskip

[H] J. Huntley.  Comparison of Maa\ss{} wave forms. {\it Number theory (New York, 1989/1990),} 129--147, Springer, New York, 1991.

\bigskip

[I-S] H. Iwaniec and P. Sarnak. Perspectives on the analytic theory of $L$-functions. GAFA 2000 (Tel Aviv, 1999).  {\it Geom. Funct. Anal.}  {\bf 2000},  Special Volume, Part II, 705--741.

\bigskip

[J-PS-S] H. Jacquet, I.I. Piatetskii-Shapiro, and J.A. Shalika. Rankin-Selberg convolutions. Amer. J. Math. {\bf 105} (1983), no. 2, 367--464.

\bigskip

[L] E. Lapid.  On the fine spectral expansion of Jacquet's relative trace formula.  Preprint available at http://www.math.huji.ac.il/~erezla/papers/expansionfin.pdf

\bigskip

[L-R-S] W. Luo, Z. Rudnick and P. Sarnak, {\it On the generalized Ramanujan conjecture for ${\rm GL}(n)$,} in Automorphic forms, automorphic representations, and arithmetic, Proc. Sympos. Pure Math., vol. 66, Part 2, Amer. Math. Soc., Providence, RI, 1999, pp. 301-310.

\bigskip

[Ma] I.G. Macdonald. Symmetric functions and Hall polynomials. Second edition. With contributions by A. Zelevinsky. Oxford Mathematical Monographs. Oxford Science Publications. {\it The Clarendon Press, Oxford University Press, New York,} 1995. x+475 pp.

\bigskip

[Molt] G. Molteni, {\it Upper and lower bounds at $s=1$ for certain Dirichlet series with Euler product}, Duke Math J., Vol. 111, No. 1 (2000), 133-158.

\bigskip

[Mo1] C. Moreno. Explicit Formulas in the theory of automorphic forms, {\it Number Theory Day} (Proc. Conf., Rockefeller Univ., New York, 1976), pp.73-216.  Lecture Notes in Math., Vol. 626, {\it Springer, Berlin,} 1977.

\bigskip

[Mo2] C. Moreno. Analytic proof of the strong multiplicity one theorem, {\it Amer. J. Math.} {\bf 107} (1985), no.1, 163-206.

\bigskip

[Mu] R. Murty. Congruences between modular forms, {\it Analtic Number Theory (Kyoto, 1996)}, 309-320, London. Math.Soc. Lecture Notes Ser., 247, Cambridge University Press, Cambridge, 1997.

\bigskip

[R-W] D. Ramakrishnan and S. Wang. On the exceptional zeros of Rankin-Selberg $L$-functions.  {\it Compositio Math.} {\bf 135} (2003), no. 2, 211--244.

\bigskip

[PS] I.I. Piatetski-Shapiro.   Multiplicity one theorems, {\it Proc. Sympos. Pure Math.}, {\bf 33}, Part 1 (1979), 209-212.

\bigskip

[R-S] Z. Rudnick and P. Sarnak.  Zeros of principal $L$-functions and random matrix theory. A celebration of John F. Nash, Jr. {\it Duke Math. J.} {\bf 81} (1996), no.2, 269-322.

\bigskip

[S] P. Sarnak. Nonvanishing of $L$-functions on $\Re (s)=1$.  {\it A supplemental volume to the Amer. J. Math.,}  Contributions to Automorphic forms, Geometry and Number Theory: Shalikafest 2002, 2003.

\end{document}